\documentclass[a4paper,11pt]{amsart}
\usepackage{amsfonts,amssymb,amsthm,cite,amsmath,amstext}
\usepackage{color}
\usepackage[colorlinks,linkcolor=black,citecolor=black]{hyperref}
\usepackage{comment}
\usepackage{framed}
\usepackage{calligra} 
\usepackage{tikz-cd}

\definecolor{shadecolor}{gray}{0.875}

\newtheorem{thrm}{Theorem}[section]
\newtheorem{thrmx}{Theorem}

\newtheorem{lem}[thrm]{Lemma}
\newtheorem{cor}[thrm]{Corollary}
\newtheorem{prop}[thrm]{Proposition}

\theoremstyle{definition}
\newtheorem{defn}[thrm]{Definition}

\newtheorem{rmk}[thrm]{Remark}

\newtheorem{ques}[thrm]{Question}

\DeclareMathOperator{\nd}{nd}

\DeclareMathOperator{\HL}{HL}
\DeclareMathOperator{\HR}{HR}

\DeclareMathOperator{\prim}{Prim}
\DeclareMathOperator{\imag}{Image}
\DeclareMathOperator{\bs}{Bs}

\title{Hard Lefschetz theorems for free line bundles}
\author{Jiajun Hu, Shijie Shang, Jian Xiao}
\date{}

\begin{document}
\maketitle


\begin{abstract}
We introduce a partial positivity notion for algebraic maps via the defect of semismallness. This positivity notion is modeled on $m$-positivity in the analytic setting and $m$-ampleness in the geometric setting. Using this positivity condition for algebraic maps, we establish K\"ahler packages, that is, Hard Lefschetz theorems and Hodge-Riemann bilinear relations, for the complete intersections of Chern classes of free line bundles.
\end{abstract}

\tableofcontents

\section{Introduction}
In this paper, we work on the field of complex numbers $\mathbb{C}$.

\subsection{Motivation}
Given a mathematical object $X$ of ``dimension'' $n$, a K\"ahler package on $X$ usually consists of a triple $(A(X), P(X), K(X))$, where $A(X)$ is a graded vector space (or even a graded algebra) constructing from $X$, $P(X)$ is a bilinear paring on $A(X)$ and $K(X)$ is a family of linear operators acting on $A(X)$ of degree one. The triple is called a K\"ahler package if it satisfies Poincar\'{e} duality, hard Lefschetz theorem and Hodge-Riemann bilinear relation. For example, given $q\leq n/2$, the hard Lefschetz theorem for $X$ means the following statement: for any $L_1, ..., L_{n-2q}\in K(X)$, the linear map
  \begin{equation*}
    A^q (X) \rightarrow A^{n-q} (X),\ \xi \mapsto \left(\prod_{k=1} ^{n-2q} L_k \right) \cdot \xi
  \end{equation*}
is an isomorphism.
In last decades, novel and exciting K\"ahler packages were discovered and have been playing key roles in algebra, combinatorics and geometry, see e.g. \cite{huhICM2018, huh2022combinatorics, williamECM} and the references therein.

As a sequel to \cite{huxiaohardlef2022Arxiv}, we are interested in the question: for which kinds of $L_1,...,L_{n-2q}$ the hard Lefschetz theorem and Hodge-Riemann bilinear relation will hold. One motivation is that in certain problems one indeed needs to study K\"ahler packages without strong positivity assumption. Though in different context, see e.g. Adiprasito's recent breakthrough on the g-conjecture of McMullen in full generality \cite{adiprasito2019combinatorial} where hard Lefschetz theorem beyond positivity is an essential ingredient.

We focus on the model and original case -- the K\"ahler package on a compact K\"ahler manifold, which could be the prototype for later development.
Let $X$ be a compact K\"ahler manifold of dimension $n$ and $\omega$ a K\"ahler class on $X$. As a fundamental piece of Hodge theory, for any integers $0\leq p, q \leq p+q \leq n$, the complete intersection class $$\Omega = \omega^{n-p-q} \in H^{n-p-q, n-p-q} (X, \mathbb{R}),$$
has the following properties:
\begin{description}
  \item [(HL)] The linear map $$\Omega: H^{p,q}(X, \mathbb{C})\rightarrow H^{n-q,n-p}(X, \mathbb{C}),\ \phi \mapsto \Omega \cdot \phi$$
      is an isomorphism.
  \item [(HR)] The quadratic form $Q$ on $H^{p, q}(X, \mathbb{C})$, defined by
  \begin{equation*}
    Q(\varphi_1, \varphi_2) = c_{p, q} \Omega \cdot  \varphi_1 \cdot \overline{\varphi_2}, \ \text{where}\ c_{p, q} = \mathrm{i}^{q-p} (-1)^{(p+q)(p+q+1)/2},
  \end{equation*}
  is positive definite on the primitive space $\prim^{p, q}(X, \mathbb{C})$ with respect to $(\Omega, \omega)$:
  \begin{equation*}
    \prim^{p, q}(X, \mathbb{C})=\{\phi \in H^{p, q}(X, \mathbb{C})| \Omega \cdot \omega\cdot \phi=0.\}
  \end{equation*}

\end{description}

We call that \emph{the class $\Omega$ has HL property} and \emph{the pair $(\Omega, \omega)$ has HR property}.

By replacing the above $\Omega=\omega^{n-p-q}$ by an arbitrary cohomology class $\Omega \in H^{n-p-q, n-p-q} (X, \mathbb{R})$ and $\omega$ by an arbitrary $(1,1)$ class $\eta \in  H^{1, 1} (X, \mathbb{R})$,
it is interesting to study:
\begin{quote}
  When does the class $\Omega$ have HL property and when does the pair $(\Omega, \eta)$ have HR property? (In order to define the primitive space we need the class $\Omega$ coupling with an $(1,1)$ class $\eta$.)
\end{quote}

A complete characterization of such $\Omega$ seems unreachable at this moment, nevertheless, in the same spirit of \cite{huxiaohardlef2022Arxiv} we first study a subclass of $\Omega$, coming from complete intersections:

\begin{ques}\label{ques HL}
Let $X$ be a compact K\"ahler manifold of dimension $n$, and let $\alpha_1,....,\alpha_{n-p-q}, \eta \in H^{1,1}(X, \mathbb{R})$, then under which assumptions does the complete intersection class $$\Omega=\alpha_1 \cdot.... \cdot\alpha_{n-p-q}$$ have HL property and does the pair $(\Omega, \eta)$ have HR property?
\end{ques}

In \cite{huxiaohardlef2022Arxiv}, under a mild positivity assumption (nefness) on the $(1,1)$ classes, the first and third named authors gave a complete characterization of $\Omega$ on a compact complex torus:

\begin{thrm} [Theorem A of \cite{huxiaohardlef2022Arxiv}]\label{hardlef tori}
Let $X=\mathbb{C}^n / \Gamma$ be a compact complex torus of dimension $n$ and $0\leq p,q\leq p+q\leq n$. Let $\alpha_1,...,\alpha_{n-p-q}, \eta\in H^{1,1} (X, \mathbb{R})$ be nef classes on $X$. Denote $\Omega=\alpha_1\cdot...\cdot\alpha_{n-p-q}$. Then for HL property, the following statements are equivalent:
\begin{enumerate}
  \item the intersection class $\Omega$ has HL property;
  \item for any subset $I \subset [n-p-q]$, $\nd(\alpha_I)\geq |I|+p+q$.
\end{enumerate}
For HR property, the following statements are equivalent:
\begin{enumerate}
  \item the pair $(\Omega, \eta)$ has HR property for $\eta$ with $\nd(\eta)\geq p+q$;
  \item for any subset $I \subset [n-p-q]$, $\nd(\alpha_I)\geq |I|+p+q$.
\end{enumerate}

\end{thrm}

Here, $[k]$ is the finite set $\{1,2,...,k\}$ for a given positive integer $k$, $\nd(-)$ is the numerical dimension of nef classes, $\alpha_I =\sum_{i\in I} \alpha_i$ and $|I|$ is the cardinality of $I$.

As a consequence of Theorem \ref{hardlef tori}, we obtain new kinds of cohomology classes on an arbitrary compact K\"ahler manifold, which have HL and HR properties.

\begin{cor}[Corollary A of \cite{huxiaohardlef2022Arxiv}]\label{cor kahler hrr intro}
Let $X$ be a compact K\"ahler manifold of dimension $n$ and $0\leq p,q\leq p+q\leq n$. Let $\alpha_1,...,\alpha_{n-p-q}, \eta \in H^{1,1}(X,\mathbb{R})$ be nef classes on $X$. Denote $\Omega=\alpha_1\cdot...\cdot\alpha_{n-p-q}$.
Assume that there exists a smooth semi-positive representative $\widehat{\alpha}_i$ in each class $\alpha_i$, such that for any subset $I \subset [n-p-q]$, $\widehat{\alpha}_I$ is $|I|+p+q$ positive in the sense of forms, and that there exists a smooth semi-positive representative $\widehat{\eta}$ in $\eta$, such that $\widehat{\eta}$ is $p+q$ positive, then
\begin{itemize}
  \item the complete intersection class $\Omega$ has HL property;
  \item the pair $(\Omega, \eta)$ has HR property.
\end{itemize}
\end{cor}

The above results greatly generalize \cite{xiaomixedHRR, DN06, cattanimixedHRR, timorinMixedHRR} by allowing degenerate positivity for each $(1,1)$ class.
A smooth semi-positive $(1,1)$ form is $m$-positive if and only if its coefficient matrix has at least $m$ positive eigenvalues. As a typical example, if $f:X\rightarrow Y$ is a submersion from $X$ to a compact K\"ahler manifold $Y$ of dimension $m$, then for any K\"ahler class $\omega_Y$ on $Y$, the pullback $f^* \omega_Y$ is $m$-positive on $X$.

From the geometric viewpoint of analytic/algebraic maps, in the above example $f$ being a submersion looks quite restrictive. Given a proper surjective holomorphic map $g:X\rightarrow Y$, even if we assume that
$$\dim g^{-1} (y) \leq {n-m}$$
for any $y\in Y$, it is not sufficient to guarantee the $m$-positivity of $g^* \omega_Y$. On the other hand, it is easy to see that the requirement on numerical dimensions in Theorem \ref{hardlef tori} is not sufficient for a general variety. 
Nevertheless, de Cataldo-Migliorini \cite{decataldoLefsemismall} proved the following hard Lefschetz property: if $L$ is a line bundle on a complex projective manifold $X$ of dimension $n$, such that a positive power of $L$ is generated by its global sections, then the class $\Omega=c_1 (L) ^{n-p-q}$
has HL property if and only if $L$ is lef, in the sense that
$mL=f^*A$
for some projective semismall morphism $f: X\rightarrow Y$, where $A$ is an ample line bundle on $Y$. In particular, the pair $(\Omega, c_1 (L))$ has HR property. The map $f: X\rightarrow Y$ is called semismall if for every $k\geq 0$,
$$\dim Y^k + 2k \leq \dim X,$$
where $Y^k = \{y\in Y | \dim f^{-1} (y) =k\}$. This result had been playing important roles in their geometric study of Hodge theory of algebraic maps \cite{decataldoLefsemismall, catalMigHodgeMaps}.

\subsection{The main result}
Inspired by \cite{decataldoLefsemismall} and a proposal proposed by the third named author \cite{xiaomixedHRR}, from the viewpoint of K\"ahler packages, we make the following analog:

\bigskip

\begin{center}

\begin{tabular}{|c|c|c|}

  forms/classes & bundles & maps \\
  \hline
  $m$-positivity & $(n-m)$-ample & $m$-lef \\

\end{tabular}
.
\end{center}

\bigskip

The notion $m$-positivity is frequently studied in geometric partial differential equations, and the $m$-ampleness of a line bundle was studied by Sommese \cite{sommesMample}. The $m$-lefness for maps is defined via the defect of semismallness:

\begin{defn}\label{intr defn mlef line bdl}
Let $X$ be an analytic variety and $L$ a free line bundle on $X$, then we call the line bundle $L$ \emph{$m$-lef} if the Kodaira map
\begin{equation*}
  \Phi_{L}: X \rightarrow Y_{L} \subset \mathbb{P}(H^0 (X, L))
\end{equation*}
is $m$-lef in the sense that the defect of semi-smallness of $ \Phi_{L}$,
  $$r(\Phi_{L}) = \max_i \{\dim Y_L ^ i + 2i - \dim X\} \leq \dim X-m$$
where $Y_L$ is the image of $X$, $Y_L ^ i = \{y\in Y | \dim f ^{-1} (x) =i \}$ and we set $\dim Y_L ^ i = -\infty$ if $Y_L^i = \emptyset$.
\end{defn}

In particular, if $m=\dim X$, then this is exactly the notion of lefness introduced by de Cataldo-Migliorini \cite{decataldoLefsemismall}. By definitions, it is easy to see: if the free line bundle $L$ is $(\dim X-m)$-ample, then it must be $m$-lef.

\begin{rmk}\label{trancendental mlef}
It is natural to extend Definition \ref{intr defn mlef line bdl} to the analytic setting as follows. Let $X$ be a compact K\"ahler manifold and $\alpha\in H^{1,1}(X, \mathbb{R})$, then we call that $\alpha$ is $m$-lef if there is a proper surjective holomorphic map $f: X\rightarrow Y$ to a K\"ahler variety $Y$ such that
\begin{itemize}
  \item $f$ is $m$-lef;
  \item $\alpha=f^* \omega_Y$ for some K\"ahler class $\omega_Y$ on $Y$.
\end{itemize}
\end{rmk}

For notional simplicity, we use the following notation: let $L_1,...,L_k$ be line bundles, then the complete intersection of their Chern classes $c_1(L_1)\cdot...\cdot c(L_k)$ is simply denoted by $L_1\cdot...\cdot L_k$.

Inspired by Corollary \ref{cor kahler hrr intro} and using the notion of $m$-lefness, we prove the following result:

\begin{thrmx}\label{intr main thrm}
Let $X$ be a smooth projective variety of dimension $n$ and let $0\leq p,q\leq p+q\leq n$ be integers. Assume that $L_1,...,L_{n-p-q}, M$ are free line bundles such that 
\begin{itemize}
  \item $L_I$ is $|I|+p+q$ lef for any $I\subset [n-p-q]$,
  \item $M$ is $p+q$ lef,
\end{itemize}
then the following statements hold:
\begin{enumerate}
  \item the complete intersection class $\Omega=L_1\cdot...\cdot L_{n-p-q}$ has HL property, i.e, the linear map
      \begin{equation*}
        \Omega: H^{p,q} (X, \mathbb{C})\rightarrow H^{n-q,n-p} (X, \mathbb{C})
      \end{equation*}
      is an isomorphism.
  \item the pair $(\Omega, M)$ has HR property, i.e., the quadratic form $Q$ on $H^{p,q} (X, \mathbb{C})$,
  \begin{equation*}
    Q(\alpha,\beta)=c_{p,q} \Omega\cdot \alpha\cdot \overline{\beta},\ \alpha, \beta\in H^{p,q} (X, \mathbb{C}),\ c_{p, q} = \mathrm{i}^{q-p} (-1)^{(p+q)(p+q+1)/2},
  \end{equation*}
  is positive definite on the primitive space
  $$\prim^{p,q} (X) =\ker \{\Omega\cdot M : H^{p,q} (X, \mathbb{C})\rightarrow H^{n-q+1,n-p+1} (X, \mathbb{C})\}.$$
\end{enumerate}
\end{thrmx}

\begin{rmk}
When $p=q=1$, by using Theorem \ref{intr main thrm} and the dually Lorentzian polynomials \cite{ross2023duallylorent}, one can obtain a bunch of generalized Alexandrov-Fenchel inequalities.
\end{rmk}

As the converse to Question \ref{ques HL}, analogous to Theorem \ref{hardlef tori}, it is interesting to study the converse of Theorem \ref{intr main thrm}: 

\begin{quote}
\emph{If the complete intersection class $\Omega=L_1\cdot...\cdot L_{n-p-q}$ has HL property, then what can we say about the positivity of the free bundles $L_i$?}
\end{quote}

Regarding this problem, we note that the HL property for certain degree is sufficient to characterize the $m$-lefness of a free line bundle, see Theorem \ref{charc mlef} for details.

We expect that Theorem \ref{intr main thrm} extends to the analytic setting as mentioned in Remark \ref{trancendental mlef}. However, the geometric arguments cannot easily extend to this analytic setting. We hope to return to this issue in the future.

The paper is organized as follows. In Section \ref{sec mlef}, we introduce the notion $m$-lefness for maps and bundles, and study its basic properties. In Section \ref{sec hironaka}, we study a variant of Hironaka's  principle of counting constants for free line bundles and the restriction of $m$-lef line bundles on a hypersurface, which will be a key ingredient in the proof of the main result. Section \ref{sec proof main theorem} is devoted to the proof of Theorem \ref{intr main thrm}.

\subsection*{Acknowledgements}
This work is supported by the National Key Research and Development Program of China (No. 2021YFA1002300) and National Natural Science Foundation of China (No. 11901336). We would like to thank Izzet Coskun and Zhiyu Tian for helpful discussions on Hironaka's original paper on the principle of counting constants, and thank Mark Andrea de Cataldo and Julius Ross for helpful comments.

\section{Partial positivity for maps and free bundles}\label{sec mlef}

The key notion in our study of partial positivity for maps is the defect of semismallness introduced by Goresky-MacPherson \cite{goreskymacphersonMorse}.

\begin{defn}
Let $f: X \rightarrow Y$ be a proper surjective holomorphic map between two complex analytic varieties, the defect of semi-smallness of $f$ is defined by
  $$r(f) = \max_i \{\dim Y ^ i + 2i - \dim X\}$$
where $Y ^ i = \{y\in Y | \dim f ^{-1} (x) =i \}$ and we set $\dim Y ^ i = -\infty$ if $Y^i = \emptyset$.
\end{defn}

The map $f$ is called semismall if $r(f)=0$.
Therefore, the number $r(f)$ measures the deviation of $f$ from semismallness.

It is clear that each $Y^i$ is a locally closed analytic subvariety of $Y$, whose disjoint union is $Y$. By letting $i=\dim X/Y$ (the dimension of a general fiber), we see that
$$r(f)\geq \dim X/Y \geq 0.$$
By the definition of defect of semi-smallness, we have that
\begin{equation}\label{eq defn defect}
  r(f)=\max_{T \subset X} \{2\dim T - \dim f(T) -\dim X\},
\end{equation}
where $T$ ranges over all irreducible analytic subvarieties (including $X$ itself) of $X$.

In the study of K\"ahler packages with respect to degenerate positivity, the third named author \cite{xiaomixedHRR} proposed an analog of $m$-positivity for analytic maps. Here, we make it more precise:

\begin{defn}
Let $f: X \rightarrow Y$ be a proper surjective holomorphic map between two complex analytic varieties and $0\leq m \leq \dim X$, then $f$ is called \emph{$m$-lef} if
$$r(f)\leq \dim X - m.$$
It is called \emph{exact $m$-lef} if $r(f)=\dim X - m$.
\end{defn}

By definition, $f$ is semismall if and only if it is $(\dim X)$-lef.

\begin{lem}\label{defect comp}
Let $f: X \rightarrow Y$, $g: X\rightarrow Z$ and $\pi: Y\rightarrow Z$ be proper surjective holomorphic maps between complex analytic varieties such that $g=\pi \circ f$, then $r(f)\leq r(g)$. Moreover, if $\pi$ is a finite morphism, then $r(f)=r(g)$.
\end{lem}

\begin{proof}
Note that for any irreducible analytic subvariety $T$ in $X$,
$$\dim g(T) =\dim \pi \circ f(T) \leq \dim f(T),$$
therefore by (\ref{eq defn defect}) we have $r(f)\leq r(g)$.

Furthermore, assume that $\pi$ is a finite morphism, then for any $i\geq 0$, $Z^i \subset \pi(Y^i)$, which yields that $r(g)\leq r(f)$. Thus, using the above inequality, $r(f)=r(g)$.
\end{proof}

Let $X$ be a projective variety and $|W|$ a linear series on $X$. We denote the Kodaira map corresponding to $W$ by
\begin{equation*}
  \Phi_W: X\dashrightarrow \mathbb{P}(W),
\end{equation*}
and denote $Y_W$ the Zariski closure of $\Phi_W (X\setminus \bs |W|)$. In particular, if the linear series is given by a line bundle $L$, we shall use the notations $\Phi_L, Y_L$.

\begin{defn}\label{defn mlef line bdl}
Let $X$ be a smooth projective variety and $L$ a free line bundle on $X$, then we call the line bundle $L$ \emph{$m$-lef} if the Kodaira map
\begin{equation*}
  \Phi_{L}: X \rightarrow Y_{L} \subset \mathbb{P}(H^0 (X, L))
\end{equation*}
is $m$-lef. The bundle $L$ is called exact $m$-lef if $\Phi_L$ is exact $m$-lef.
\end{defn}

When $m=\dim X$, this is exactly the notion \emph{lefness} introduced by de Cataldo-Migliorini \cite{decataldoLefsemismall} in their study of decomposition theorem for semismall maps.

\begin{rmk}
By \cite{viehwegvanishingnonvan}, Kodaira-Akizuki-Nakano type vanishing theorems hold with respect to $m$-lefness: let $X$ be a smooth projective variety and let $L$ be $m$-lef on $X$, then
\begin{equation*}
  H^{p,q} (X, L^{-1}) = 0
\end{equation*}
whenever $p+q< m$.
\end{rmk}

\begin{prop}\label{num dim lef}
Assume that the line bundle $L$ is $m$-lef, then its numerical dimension $\nd(L)\geq m$.
\end{prop}

\begin{proof}
The follows by letting $i=\dim X/ Y_L$ (the dimension of a general fiber of $\Phi_L$) in the formula of $r(\Phi_L)$. Then we get $\dim Y_L \geq m$, which means that the Kodaira dimension of $L$ is at least $m$.

For a free line bundle, note that its Kodaira dimension and numerical dimension coincide.
\end{proof}

\begin{rmk}
By Proposition \ref{num dim lef}, the positivity assumption in Theorem \ref{intr main thrm} implies the requirement on numerical dimensions as in Theorem \ref{hardlef tori}.
\end{rmk}

The following result is useful in the comparison of partial lefness for different free linear series.

\begin{lem}\label{finite proj}
Let $X$ be a smooth projective variety, and $|W|, |V|$ two free linear series on $X$. Assume that $|W| \subset |V|$, then the two morphisms $\Phi_W: X\rightarrow \mathbb{P}(W)$, $\Phi_V: X\rightarrow \mathbb{P}(V)$ differ by a finite projection of $\Phi_V (X)$, that is, $\Phi_W =\pi \circ \Phi_V$ where $\pi$ is a finite morphism.
\end{lem}

\begin{proof}
This follows from \cite[Example 1.1.12]{lazarsfeldPosI}.
\end{proof}

For example, let $L$ be a free line bundle, by Lemma \ref{finite proj}, it is clear that the following are equivalent:
\begin{itemize}
  \item $L$ is $m$-lef;
  \item $kL$ is $m$-lef for some positive integer $k$;
  \item $kL$ is $m$-lef for any positive integer $k$.
\end{itemize}

\begin{prop}
Any nontrivial free line bundle $L$ is $m$-lef for some $m\geq 1$.
\end{prop}

\begin{proof}

Consider the Kodaira map
$$\Phi_L: X\rightarrow Y_L \subset \mathbb{P}(H^0 (X, L)).$$
It is easy to see that
\begin{equation*}
  2\dim T - \dim \Phi_L (T) \leq 2\dim X - 1
\end{equation*}
holds for any irreducible subvariety $T$ of $X$. Therefore, any nontrivial free line bundle is $1$-lef. Indeed, $L$ is exact $m$-lef, where
\begin{equation*}
 m=\min_{T\subset X} \{2\dim X-2\dim T + \dim \Phi_L (T)\}.
\end{equation*}

\end{proof}

Therefore, the notion $m$-lefness gives a filtration on the positivity of free line bundles.

\begin{lem}\label{sum lef max}
If $K$ is $k$-lef and $L$ is $l$-lef, then $K+L$ is $\max(k, l)$-lef.
\end{lem}

\begin{proof}
By comparing the linear series $|K+L|$ and $|K|\otimes |L|$ and using Lemma \ref{finite proj}, it suffices to show that the map
\begin{equation*}
  f: X\rightarrow Y_K \times Y_L,\ x\mapsto (\Phi_K (x), \Phi_L (x)),
\end{equation*}
is $\max(k, l)$-lef onto its image.
Applying Lemma \ref{defect comp}, we have that
$$r(f)\leq \min\{r(\Phi_K), r(\Phi_L)\}.$$
This yields that $f$ is $\max(k, l)$-lef.

\end{proof}

The following Lefschetz hyperplane theorem due to Goresky-MacPherson \cite[Section 2.3]{goreskymacphersonMorse} is useful for us.
\begin{lem}\label{weakLef}
Let $X$ be a smooth projective variety and let $L$ be a $r$-lef free line bundle on $X$. Then for a general hypersurface $V\in |L|$, the restriction map
\begin{equation*}
  i^*: H^{l}(X,\mathbb{C})\rightarrow H^{l}(V,\mathbb{C})
\end{equation*}
is injective for $l \leq r-1$ and isomorphic for $l \leq r-2$.
\end{lem}

As a direct application, we get the following Bertini theorem for $2$-lef line bundles.

\begin{lem}\label{BertiniLef}
  Let $X$ be a smooth projective variety and let $L$ be a $2$-lef free line bundle on $X$. Then a general hypersurface $V\in |L|$ is smooth and irreducible.
\end{lem}

\begin{rmk}
We give a brief discussion on the background for the analog:

\begin{center}
\bigskip
\begin{tabular}{|c|c|c|}

  forms/classes & bundles & maps \\
  \hline
  $m$-positivity & $(n-m)$-ample & $m$-lef \\

\end{tabular}
,
\bigskip
\end{center}
and the motivation for $m$-lefness.
Let $\omega$ be a K\"ahler class on a compact K\"ahler manifold $X$ of dimension $n$ and let $\widehat{\omega}$ be a K\"ahler metric in the class $\omega$. Then $\alpha\in H^{1,1} (X, \mathbb{R})$ is called $m$-positive with respect to $\widehat{\omega}$, if $\alpha$ has a smooth representative $\widehat{\alpha}$ such that for any $1\leq k\leq m$,
\begin{equation}\label{eq mposi}
  \widehat{\alpha}^k \wedge \widehat{\omega}^{n-k}>0
\end{equation}
in the sense of forms. If we further assume that $\widehat{\alpha}$ is semipositive, (\ref{eq mposi}) is equivalent to that $\widehat{\alpha}$ has at least $m$ positive eigenvalues everywhere with respect to $\widehat{\omega}$. In particular, a free line bundle $L$ is called $m$-positive with respect to $\widehat{\omega}$, if its Chern class $c_1 (L)$ is $m$-positive with respect to $\widehat{\omega}$. In \cite{xiaomixedHRR}, the third named author proved the following result:
\begin{quote}
  Let $0\leq p,q\leq p+q\leq m\leq n$ be integers, and let $\alpha_1,...,\alpha_{m-p-q+1}\in H^{1,1} (X, \mathbb{R})$ be semipositive and $m$-positive classes, then the class $\Omega=\omega^{n-m}\cdot\alpha_1\cdot...\cdot\alpha_{m-p-q}$ has HL property and the pair $(\Omega, \alpha_{m-p-q+1})$ has HR property.
\end{quote}
Applying the result to the case when every $\alpha_k$ coming from a submersion $f:X\rightarrow Y$ motivates essentially the following notion of partial lefness for line bundles \cite[Section 4]{xiaomixedHRR}:
a free line bundle $L$ is called $m$-lef if the Kodaira map
\begin{equation*}
  \Phi_{L}: X \rightarrow Y_{L}
\end{equation*}
satisfies that
\begin{itemize}
  \item the numerical dimension of $L$, $\nd(L) \geq m$;
  \item the defect of semi-smallness of $f$,
  $r(f)\leq n - \dim Y_L $.

\end{itemize}
In general, this notion is stronger than Definition \ref{defn mlef line bdl}.

In \cite{sommesMample}, Sommese introduced the notion $m$-ampleness for line bundles.
Let $X$ be a smooth projective variety of dimension $n$ and $L$ a line bundle on $X$, the line bundle $L$ is called $m$-ample if there exists some $k\in \mathbb{N}$ such that $\dim \bs_{|kL|} \leq m$ and the Kodaira map
\begin{equation*}
  \Phi_{kL}: X \setminus \bs_{|kL|} \rightarrow Y_{kL}\subset \mathbb{P}(H^0 (X, kL))
\end{equation*}
satisfies that $\dim \Phi_{kL} ^{-1} (x) \leq m$ for any $x$ in the image.
In particular, if $L$ is free, then $L$ is $m$-ample if and only if $\dim \Phi_{L} ^{-1} (x) \leq m$ for any $x\in Y_L$.

In the case when $\Phi_{kL}$ is a submersion, we have:
\begin{center}
$L$ being $(n-m)$-ample $\Rightarrow$ $c_1 (L)$ being semipositive and $m$-positive.
\end{center}
By Definition \ref{defn mlef line bdl}, for free bundles we also have:
\begin{center}
$L$ being $(n-m)$-ample $\Rightarrow$ $L$ being $m$-lef.
\end{center}

\end{rmk}

Next, similar to \cite{decataldoLefsemismall} we give a characterization of $m$-lefness by the hard Lefschetz property.

\begin{prop}
Let $X$ be a smooth projective variety of dimension $n$ and $L$ a free line bundle on $X$. Assume that $L$ is $m$-lef, then
for any ample line bundles $A_1,...,A_{n-m}$, the class $A_1\cdot ...\cdot A_{n-m}\cdot L^{m-p-q}$ has HL property for any $0\leq p, q\leq p+q\leq m$, that is, the map
\begin{equation*}
  A_1\cdot...\cdot A_{n-m}\cdot L^{m-p-q}: H^{p,q}(X, \mathbb{C})\rightarrow H^{n-q,n-p}(X, \mathbb{C})
\end{equation*}
is an isomorphism.
\end{prop}

\begin{proof}
For $p+q =m$, the result is clear, thus we need only to deal with the cases $p+q\leq m-1$. After taking multiples of the line bundles we can assume that $A_1,...,A_{n-m}$ are very ample.

Take a general smooth hypersurface $H_i \in |A_i|$ for each $i\leq n-m$, then by \cite[Section 4]{catalMigHodgeMaps}, the restriction $\Phi_L$ on $V=H_1 \cap...\cap H_{n-m}$ has vanishing defect of semismallness, in particular it must be semismall. This implies that the restriction of $L$ on $V$, denoted by $L_{|V}$, is lef.

Assume that $\varphi \in H^{p+q} (X, \mathbb{C})$ satisfies
\begin{equation*}
  A_1\cdot ...\cdot A_{n-m}\cdot L^{m-p-q} \cdot \varphi =0,
\end{equation*}
we need to prove that $\varphi=0$. To this end, note that this is equivalent to that
\begin{equation}\label{restr}
  A_1\cdot ...\cdot A_{n-m}\cdot L^{m-p-q} \cdot \varphi \cdot \psi = L_{|V} ^{m-p-q} \cdot \varphi_{|V} \cdot \psi_{|V} =0
\end{equation}
for any $\psi \in H^{p+q} (X, \mathbb{C})$. By an inductive application of the Lefschetz hyperplane theorem, we obtain that the map
\begin{equation*}
  i_V: H^{p+q} (X, \mathbb{C}) \rightarrow H^{p+q} (V, \mathbb{C}),\ \psi\mapsto \psi_{|V},
\end{equation*}
is an isomorphism whenever $p+q \leq m-1$. Combining with (\ref{restr}), we get that
\begin{equation*}
  L_{|V} ^{m-p-q} \cdot \varphi_{|V} =0
\end{equation*}
on $V$.

By \cite{decataldoLefsemismall}, since $L_{|V}$ is lef on $V$, it has HL property. Therefore, $\varphi_{|V}=0$. Another application of the Lefschetz hyperplane theorem implies that $\varphi =0$, finishing the proof.

\end{proof}

\begin{prop}
Let $X$ be a smooth projective variety and $L$ a free line bundle on $X$. Let $A_1,...,A_{n-m}$ be ample line bundles on $X$. Assume that the class $A_1\cdot ...\cdot A_{n-m}\cdot L^{m-p-q}$ has HL property for any $0\leq p, q\leq p+q\leq m$, then $L$ is $m$-lef.
\end{prop}

\begin{proof}
Since $L$ is free, we may assume that $L= \Phi_L ^* A$ for some ample line bundle $A$ on $Y_L$.

We argue by contradiction. Otherwise, $L$ is not $m$-lef, which means that there is some irreducible subvariety $T$ in $X$ such that
\begin{equation*}
  2\dim T -2n+m> \dim f(T).
\end{equation*}
This implies that there exists a cycle in the class $A^{2\dim T -2n+m}$, which is disjoint with $f(T)$. Thus, there is a cycle in the class $$L^{2\dim T -2n+m}=\Phi_L ^* A^{2\dim T -2n+m}$$
 which is disjoint with $T$, yielding that
$$A_1\cdot ...\cdot A_{n-m}\cdot L^{m-2(n-\dim T)} \cdot [T] =0.$$
Therefore, $A_1\cdot ...\cdot A_{n-m}\cdot L^{m-p-q}$ does not have HL property for $p+q=2(n-\dim T)$.

This finishes the proof.
\end{proof}

In summary, we obtain the following characterization:

\begin{thrm}\label{charc mlef}
Let $X$ be a smooth projective variety of dimension $n$ and $L$ a free line bundle on $X$, then the following statements are equivalent:
\begin{itemize}
  \item $L$ is $m$-lef;
  \item for any ample line bundles $A_1,...,A_{n-m}$, the class $A_1\cdot ...\cdot A_{n-m}\cdot L^{m-p-q}$ has HL property for any $0\leq p, q\leq p+q\leq m$;
  \item for some ample line bundles $A_1,...,A_{n-m}$, the class $A_1\cdot ...\cdot A_{n-m}\cdot L^{m-p-q}$ has HL property for any $0\leq p, q\leq p+q\leq m$.
\end{itemize}

\end{thrm}

We expect that Theorem \ref{charc mlef} also holds in the analytic setting (see Remark \ref{trancendental mlef}).

\section{Hironaka's principle of counting constants}\label{sec hironaka}

In \cite[Section 2]{hironakaprinciple}, using the ``counting constants" method, Hironaka proved results of the following form. We refer the reader to \cite[Chapter 3]{sommeseShiffmVanishingBook} for a modern account.

\begin{thrm}[Theorem 3.39 of \cite{sommeseShiffmVanishingBook}] \label{sommese hiron}
Let $X$ be a smooth projective variety and let $L$ be an ample line bundle on $X$. Let $f:X\dashrightarrow \mathbb{P}^N$ be a meromorphic map and let $X_0$ be the Zariski open set of $X$ such that $f_0=f_{|X_0}$ is holomorphic. Then there is a positive integer $k$ such that the generic element of $H^0 (X, kL)$ does not vanish on any positive dimensional component of $f_0^{-1}(y)$ for all $y\in f_0 (X_0)$.
\end{thrm}

In our setting, we need its extension to free bundles under certain assumption on the interaction between the map and the bundle.

\begin{lem}\label{lemma1}
Let $X$ be a smooth projective variety of dimension $n$, and let $L, F$ be free line bundles on $X$ with corresponding Kodaira maps $\Phi_L, \Phi_F$.  Then there exists $m_0 \in \mathbb{N}$ such that for any $m\geq m_0$, a generic section of $mF$ does not vanish on any positive dimensional irreducible component $W$ of any fiber $\Phi_L ^{-1}(y)$ with $\dim \Phi_F (W) >0$.

\end{lem}

\begin{proof}
Let $W$ be a positive dimensional irreducible component $W$ of $\Phi_L ^{-1}(y)$ with $\dim \Phi_F (W) >0$, and
denote $q:=\dim \Phi_{F}(W)\geq 1$. We may assume that $\Phi_L:X\rightarrow Y_L$ is the Iitaka fibration associated to $L$. The proof will be divided into three steps.

\textbf{Step 1.} We have the exact sequence
\begin{equation*}
0\rightarrow H^0(X,\mathcal{I}_{W/X}\otimes mF)\rightarrow H^0(X,mF)\rightarrow H^0(X|W,mF),
\end{equation*}
where $\mathcal{I}_{W/X}$ is the ideal sheaf of $W$ in $X$ and
$$H^0(X|W,mF)=\imag [H^0(X,mF)\rightarrow H^0(W,mF|_W)]$$
is the space of restricted sections of $mF$ on $W$.

Following the proof of Pacienza-Takayama \cite[Theorem 1.1]{takayamaIntermKodai}, one can show that there exists an integer $m_0$ and a constant $c>0$ such that
$$h^0(X|W,mF)\geq cm^q$$
holds for any $m\geq m_0$.

For the convenience of readers, we include the details here. Let $A$ be the very ample line bundle on $Y_F$ such that $\Phi^*_{F}A=F$. Then we have
\begin{equation*}
H^0(Y_F,mA)\cong H^0(X,m\Phi^*_{F}A)=H^0(X,mF)
\end{equation*}
for any $m>0$. Since $A$ is ample on $Y_F$, then there exists an integer $m_0$ such that for any $m\geq m_0$ and any $i>0$, we have
\begin{equation*}
H^i(Y_F,\mathcal{I}_{\Phi_{F}(W)/Y_F}\otimes mA)=0
\end{equation*}
where $\mathcal{I}_{\Phi_{F}(W)/Y_F}$ is the ideal sheaf of $\Phi_{F}(W)$ in $Y_F$. Thus, the restriction map
\begin{equation*}
H^0(Y_F,mA)\rightarrow H^0(\Phi_{F}(W),mA)
\end{equation*}
is surjective for any $m\geq m_0$. Then we have an inclusion
\begin{equation*}
(\Phi_{F}|_W)^*H^0(\Phi_{F}(W),mA)\subset H^0(X|W,mF)
\end{equation*}
for any $m\geq m_0$. Hence there exists a constant $c>0$ such that
$$h^0(X|W,mF)\geq cm^q$$
 for any $m\geq m_0$.
~\\

The following arguments are inspired by \cite[Section 2]{hironakaprinciple}, \cite[Chapter 3]{sommeseShiffmVanishingBook} and \cite[\href{https://stacks.math.columbia.edu/tag/055A}{Tag 055A}]{Stacks}. 

\textbf{Step 2.} Moreover, $c$ and $m_0$ can be chosen independently of $W$ and $y\in Y_L$. This statement will be proved by induction on $\dim Y_L$. If $\dim Y_L=0$ holds, then the statement is trivial. Now we assume $\dim Y_L>0$.

\textbf{Claim.} We can find a morphism $v:Y\rightarrow Y_L$ with the following properties:\\
(i) $v$ is a finite open morphism,\\
(ii) $Y$ is an integral affine scheme,\\
(iii) $X\times_{Y_L}Y=\cup^N_{i=1} X_i$ is a decomposition of $X\times_{Y_L}Y$, where the fibers of the morphism $X_i\rightarrow Y$ are all geometrically integral for any $i\in\{1,2,\cdots,N\}$.  

\textbf{Proof of the Claim.} Consider the morphism $\Phi_L:X\rightarrow Y_L$, and it follows from \cite[\href{https://stacks.math.columbia.edu/tag/0551}{Tag 0551}]{Stacks} that we have the following diagram
\[
\begin{tikzcd}
X' \arrow[r, "g'"] \arrow[d, "\Phi'_L"'] & X_V \arrow[r] \arrow[d] & X \arrow[d, "\Phi_L"] \\
Y'_L \arrow[r, "g"]                      & V \arrow[r]             & Y_L                  
\end{tikzcd}
\] 
where\\
(i) $V$ is a nonempty open of $Y_L$,\\
(ii) $X_V=V\times_{Y_L} X$ and $X'=Y'_L\times_V X_V$,\\
(iii) both $g$ and $g'$ are surjective finite \'{e}tale,\\
(iv) $Y'_L$ is an irreducible affine scheme,\\
(v) all irreducible components of the generic fiber of $\Phi'_L$ are geometrically irreducible.

Denote by $\eta$ the generic point of $Y'_L$. Then all irreducible components of the generic fiber $X'_\eta$ are geometrically irreducible over $\kappa(\eta)$. Suppose that $X'_\eta=\cup^N_{i=1}X'_{i,\eta}$ is the decomposition of the generic fiber into (geometrically) irreducible components. For each $1\leq i\leq N$, let $X'_i$ be the closure of $X'_{i,\eta}$ in $X'$ and endow $X'_i$ with the induced reduced scheme structure. Note that the generic fiber of $X'_i$ is $X'_{i,\eta}$. After shrinking $Y'_L$ we may assume that $X'=\cup^{N}_{i=1}X'_i$ by \cite[\href{https://stacks.math.columbia.edu/tag/054Y}{Tag 054Y}]{Stacks}. After shrinking $Y'_L$ some more, it follows from \cite[\href{https://stacks.math.columbia.edu/tag/0554}{Tag 0554}]{Stacks} and \cite[\href{https://stacks.math.columbia.edu/tag/0559}{Tag 0559}]{Stacks} that $X'_{i,y}$ is geometrically irreducible for each $i$ and all $y\in Y'_L$, and $X'_y=\cup^{N}_{i=1}X'_{i,y}$ is the decomposition of the fiber $X'_y$ into (geometrically) irreducible components for all $y\in Y'_L$.

Fix any $i\in\{1,2,\dotsc,N\}$, we consider the morphism $\Phi'_L|_{X'_i}:X'_i\rightarrow Y_L$. By abuse of notation, we will write  $\Phi'_L$ instead of $\Phi'_L|_{X'_i}$. It follows from \cite[\href{https://stacks.math.columbia.edu/tag/0550}{Tag 0550}]{Stacks} that we have the following diagram
\[
\begin{tikzcd}
X''_i \arrow[r, "h'"] \arrow[d, "\Phi''_L"'] & {X'_{i,U}} \arrow[r] \arrow[d] & X'_i \arrow[d, "\Phi'_L"] \\
Y''_L \arrow[r, "h"]                         & U \arrow[r]                    & Y'_L                     
\end{tikzcd}
\]
where\\
(i) $U$ is a nonempty open of $Y'_L$,\\
(ii) $X'_{i,U}=U\times_{Y'_L} X$ and $X''_i=(Y''_L\times_U X'_{i,U})_{red}$,\\
(iii) both $h$ and $h'$ are finite universal homeomorphisms,\\
(iv) $Y''_L$ is an integral affine scheme,\\
(v) $\Phi''_L$ is flat and of finite presentation,\\
(vi) the generic fiber of $\Phi''_L$ is geometrically reduced.

After shrinking $Y''_L$, we can assume that for every point $y\in Y''_L$, the fiber $X''_{i,y}$ is geometrically integral by \cite[Theorem 12.2.1]{EGAiv}. Then the Claim follows immediately.
~\\

We will still use the notations introduced in the proof of the Claim. We denote by $\tilde{\Phi}_F$ the composition of the following morphisms 
\[
\begin{tikzcd}
X''_i \arrow[r, "h'"] & {X'_{i,U}} \arrow[r] & X'_i \arrow[r, hook] & X' \arrow[r, "g'"] & X_V \arrow[r] & X \arrow[r, "\Phi_F"] & Y_F.
\end{tikzcd}
\]
Since projective morphisms are preserved by base change, we can see that $\Phi''_L:X''_i\rightarrow Y''_L$ is projective. Then it follows immediately that the morphism $\tilde{\Phi}_{F}\times\Phi''_{L}:X''_i\rightarrow Y_F\times Y''_L$ is projective. In particular, the morphism $\tilde{\Phi}_{F}\times\Phi''_{L}$ is closed. Hence, $(\tilde{\Phi}_{F}\times\Phi''_{L})(X''_i)$ is a closed subset in $Y_F\times Y''_L$.

For any $i\in\{1,2,\cdots,N\}$, set $\Pi_i:=(\tilde{\Phi}_{F}\times\Phi''_{L})(X''_i)$ with the induced reduced subscheme structure in $Y_F\times Y''_L$, which is exactly the scheme-theoretic image of the morphism $\tilde{\Phi}_{F}\times\Phi''_{L}$. 
Let $q_i:\Pi_i\rightarrow \Phi''_{L}(X''_i)=Y''_L$ 
be the projection for each $i$. One can easily see that $q_i$ is a projective morphism for each $i$. It follows from the Claim that there exists a open subvariety $\tilde{Y}\subseteq Y_L$ such that for any closed point $y\in \tilde{Y}$, ${\Phi}_{F}(W)$ can be identified with $q^{-1}_i(y'')$ for some $i$ where $y''\in Y''_L$. Then the set of such $\Phi_F(W)$ is contained in a finite number of components of the Hilbert scheme on $\mathbb{P}({H^0(X,F)})$ by flattening stratification theorem (see e.g. \cite[Lecture 8]{Mumford}) and \cite[Theorem III.9.9]{HartsAG}. It follows from \cite[Chapter IX, Lemma 4.1]{GeometryOfAlgebraicCurvesVol2} and the following commutative diagram (all rows and columns are exact)
\[
\begin{tikzcd}
& 0 \arrow[d]           & 0 \arrow[d]           &                       &   \\
0 \arrow[r] & \mathcal{I}_{Y_F/\mathbb{P}({H^0(X,F)})} \arrow[d] \arrow[r] & \mathcal{I}_{Y_F/\mathbb{P}({H^0(X,F)})} \arrow[d] \arrow[r] & 0 \arrow[d]           &   \\
0 \arrow[r] & \mathcal{I}_{\Phi_{F}(W)/\mathbb{P}({H^0(X,F)})} \arrow[d] \arrow[r] & \mathcal{O}_{\mathbb{P}({H^0(X,F)})}  \arrow[d] \arrow[r] & \mathcal{O}_{\Phi_{F}(W)}  \arrow[r] \arrow[d] & 0 \\
0 \arrow[r] & \mathcal{I}_{\Phi_{F}(W)/Y_F} \arrow[d] \arrow[r] & \mathcal{O}_{Y_F} \arrow[d] \arrow[r] & \mathcal{O}_{\Phi_{F}(W)} \arrow[r] \arrow[d] & 0 \\
& 0                     & 0                     & 0                     &
\end{tikzcd}
\]
that the Castelnuovo-Mumford regularity of $\mathcal{I}_{\Phi_{F}(W)/Y_F}$ with respect to the very ample line bundle $A$ only depends on $h^0(X,F)$ and the Hilbert polynomial of $\mathcal{I}_{\Phi_{F}(W)/\mathbb{P}({H^0(X,F)})}$. 

By induction on the dimension of $Y_L$, we can conclude that $c$ and $m_0$ can be chosen independently of $W$ and $y\in Y_L$.

\textbf{Step 3.} Now we can choose a uniform $m_0$ such that for any $m\geq m_0$,
\[
\begin{split}
h^0(X,\mathcal{I}_{W/X}\otimes mF)&=h^0(X,mF)-h^0(X|W,mF)\\
&\leq h^0(X,mF)-h^0(X,L)-1.
\end{split}
\]

Consider the set $\Sigma$ of pairs $$(y,s)\in \mathbb{P}(H^0(X,L))\times \mathbb{P}(H^0(X,mF))$$ such that $s$ vanishes on some positive dimensional component of $\Phi^{-1}_{L}(y)$. Then $\Sigma$ is a projective variety and $$\dim \Sigma\leq h^0(X,mF)-2,.$$

Therefore, any section $s\in \mathbb{P}(H^0(X,mF))\backslash p(\Sigma)$ does not vanish identically on any
positive dimensional irreducible component $W$ of any fiber of $\Phi_{L} ^{-1}(y)$ with $\dim \Phi_F (W)>0$, where $$p:Y_L\times \mathbb{P}(H^0(X,mF))\rightarrow \mathbb{P}(H^0(X,mF))$$ is the projection.

\end{proof}

As an application of Lemma \ref{lemma1}, we obtain the following result on the restriction of $m$-lef line bundles.

\begin{prop}\label{restriction prop}
Let $L, F$ be free line bundles. Assume that $F$ is $2$-lef, $L$ is $k$-lef and $L+F$ is $(k+1)$-lef, then there exists $m_0 \in \mathbb{N}$ such that for any $m\geq m_0$ and $V\in |mF|$ generic, the restriction $L_{|V}$ is $k$-lef on the hypersurface $V$.
\end{prop}

\begin{proof}
Let $\Phi_{L}:X\rightarrow Y_L\subseteq\mathbb{P}({H^0(X,L)})$, $\Phi_{F}:X\rightarrow Y_F\subseteq\mathbb{P}({H^0(X,F)})$ and $$\Phi_{L+F}:X\rightarrow Y_{L+F}\subseteq\mathbb{P}({H^0(X,L+F)})$$ be the morphisms induced by the linear systems $|L|$, $|F|$ and $|L+F|$ respectively.

For any $i\geq 0$, set
\begin{equation*}
Y^i:=\{y\in Y_L|\dim\Phi^{-1}_{L}(y)=i\}
\end{equation*}
and
\begin{equation*}
Z^i:=\{(y,z)\in (\Phi_{L}\times\Phi_{F})(X)|\dim \Phi^{-1}_{L}(y)\cap\Phi^{-1}_{F}(z)=i\}\subseteq Y_L\times Y_F
\end{equation*}
Since $L$ is $k$-lef, for any $i\geq 0$ we have
\begin{equation*}
\dim Y^i\leq 2n-2i-k.
\end{equation*}
Since $\Phi_{L+F}$ and $\Phi_{F} \times \Phi_{L}$ only differ by a finite morphism and $L+F$ is $(k+1)$-lef, then for any $i\geq 0$ we have
\begin{equation*}
\dim Z^i\leq 2n-2i-k-1.
\end{equation*}
For any $i\geq 0$, set
\begin{equation*}
Y^{i'}:=\{y\in Y^i|\exists\,\textrm{a top dimensional irreducible component $W\subseteq\Phi^{-1}_{L}(y)$ such that}\,\dim\Phi_{F}(W)=0\}
\end{equation*}
and
\begin{equation*}
Y^{i''}:=\{y\in Y^i|\dim\Phi_{F}(W)>0\,\textrm{holds for all top dimensional irreducible components $W\subseteq\Phi^{-1}_{L}(y)$}\}.
\end{equation*}
Then for any $i\geq 0$ we have
\begin{equation*}
Y^i=Y^{i'}\sqcup Y^{i''}.
\end{equation*}

It follows from Lemma \ref{lemma1} that there exists a positive integer $m_0$ such that general element $V\in |mF|$ has proper intersection with any positive dimensional irreducible component $W$ of any fiber $\Phi^{-1}_{L}(y)$ with $\dim\Phi_{F}(W)>0$. Set $V_L:=\Phi_{L}(V)$ and
\begin{equation*}
V^i_L:=\{y\in V_L|\dim(\Phi_{L}|_V)^{-1}(y)=i\}.
\end{equation*}
By the choice of $V$, we have
\begin{equation}\label{lemma3(1)}
Y^{i''}\cap V^i_L=\emptyset
\end{equation}
and
\begin{equation*}
Y^{i+1''}\cap V_L\subseteq V^i_L.
\end{equation*}
Then we have
\begin{equation}\label{lemma3(2)}
\dim Y^{i+1''}\cap V_L\leq \dim Y^{i+1}\leq 2n-2-2i-k=2\dim V-2i-k.
\end{equation}
Let $p:Y_L\times Y_F\rightarrow Y_L$ be the projection. It follows from the definitions of $Y^{i'}$ and $Z^i$ that for any $i\geq 0$, the restriction of the projection
\begin{equation*}
p|_{Z^i\cap p^{-1}(Y^{i'})}:Z^i\cap p^{-1}(Y^{i'})\rightarrow Y^{i'}
\end{equation*}
is a surjective morphism. Hence, for any $i\geq 0$ we have
\begin{equation}\label{lemma3(5.5)}
\dim Y^{i'}\leq\dim Z^i\leq 2n-2i-k-1=2\dim V-2i-k+1.
\end{equation}
Then we have
\begin{equation}\label{lemma3(3)}
\dim Y^{i+1'}\cap V_L\leq \dim Z^{i+1}\leq 2n-2i-k-3\leq 2\dim V-2i-k.
\end{equation}

The set of divisors containing at least one irreducible component $W$ of $\Phi^{-1}_{L}(Y^{i'})$ for some $i$ is a finite union of linear proper subspaces of $|mF|$. Then we can assume that the general element $V$ satisfies
\begin{equation*}
\dim W\cap V<\dim W
\end{equation*}
for any irreducible component $W$ of $\Phi^{-1}_{L}(Y^{i'})$ and any $i\geq 0$. Meanwhile, we have
\begin{equation*}
\dim  \Phi_{L}(W\cap V)\leq \dim Y^{i'}
\end{equation*}
for any irreducible component $W$ of $\Phi^{-1}_{L}(Y^{i'})$ and any $i\geq 0$.

If $\dim \Phi_{L}(W\cap V)<\dim Y^{i'}$ holds for any irreducible component $W$ of $\Phi^{-1}_{L}(Y^{i'})$, then by
\begin{align*}
  Y^{i'}\cap V_L &= \Phi_L |_V ((\Phi_L |_ V)^{-1} (Y^{i'}\cap V_L)) \\
  &=\Phi_L |_V (\Phi_L ^{-1} (Y^{i'})\cap V)\\
  &=\bigcup_{W \subset \Phi_L^{-1} (Y^{i'})} \Phi_L (W\cap V),
\end{align*}
it follows from (\ref{lemma3(5.5)}) that
\begin{equation}\label{lemma3(4)}
\dim Y^{i'}\cap V_L\leq  2\dim V-2i-k.
\end{equation}

If $\dim \Phi_{L}(W\cap V)=\dim Y^{i'}$ holds for some irreducible component $W$ of $\Phi^{-1}_{L}(Y^{i'})$, then the dimension of general fiber of the morphism
\begin{equation*}
\Phi_{L}: W\cap V\rightarrow \Phi_{L}(W)\subseteq Y^{i'}\cap V_L
\end{equation*}
is at most $i-1$. This holds since $\dim \Phi_{L}(W\cap V)=\dim Y^{i'}$ implies that $\dim \Phi_{L}(W\cap V)=\dim \Phi_{L}(W)$, and the dimension of a general fiber of $\Phi_L : W \rightarrow \Phi_L (W)$ is at most $i$ by the definition of $Y^{i'}\subset Y^i$, and $\dim W\cap V =\dim W -1$ by the choice of $V$.
Hence, we have
\begin{equation}\label{lemma3(5)}
\dim Y^{i'}\cap V^i_L \leq \dim Y^{i'}-1\leq  2\dim V-2i-k.
\end{equation}

Combining (\ref{lemma3(1)}), (\ref{lemma3(2)}), (\ref{lemma3(3)}), (\ref{lemma3(4)}) and (\ref{lemma3(5)}), we can conclude that
\begin{equation*}
\dim V^i_L\leq 2\dim V-2i-k
\end{equation*}
holds for any $i\geq 0$.

Since $\Phi_{L}|_V$ and $\Phi_{L_{|V}}$ only differ by a finite morphism, by Lemma \ref{finite proj} and Lemma \ref{defect comp}, the proposition follows.
\end{proof}

\section{Proof of the main theorem}\label{sec proof main theorem}

Recall that we are going to prove:

\begin{thrm}\label{hardlefthrm}
Let $X$ be a smooth projective variety of dimension $n$ and let $0\leq p,q\leq p+q\leq n$ be integers. Assume that $L_1,...,L_{n-p-q},M$ are free line bundles such that $M$ is $p+q$ lef and $L_I$ is $|I|+p+q$ lef for any $I\subset [n-p-q]$, then the following statements hold:
\begin{enumerate}
  \item the complete intersection class $\Omega=L_1\cdot...\cdot L_{n-p-q}$ has hard Lefschetz property, i.e, the linear map
      \begin{equation*}
        \Omega: H^{p,q} (X, \mathbb{C})\rightarrow H^{n-q,n-p} (X, \mathbb{C})
      \end{equation*}
      is an isomorphism.
  \item the pair $(\Omega, M)$ has Hodge-Riemann property, i.e., the quadratic form $Q$ on $H^{p,q} (X, \mathbb{C})$,
  \begin{equation*}
    Q(\alpha,\beta)=c_{p,q} \Omega\cdot \alpha\cdot \overline{\beta},\ \alpha, \beta\in H^{p,q} (X, \mathbb{C}),
  \end{equation*}
  is positive definite on the primitive space
  $$\prim^{p,q} (X) =\ker \{\Omega\cdot M : H^{p,q} (X, \mathbb{C})\rightarrow H^{n-q+1,n-p+1} (X, \mathbb{C})\}.$$
\end{enumerate}
\end{thrm}

\begin{proof}

For $p=q=0$, we need to show the following statement: 
\begin{equation*}
  L_1\cdot...\cdot L_n>0.
\end{equation*}
Note that by Proposition \ref{num dim lef}, $L$ being $r$-lef implies that $\nd(L)\geq r$.
Then the above statement is a consequence of the positivity criterion for the intersection of nef classes \cite{huxiaohardlef2022Arxiv}:
\begin{center}
  $L_1\cdot...\cdot L_n>0$ if and only if $\nd(L_I)\geq |I|, \forall I\subset [n]$.
\end{center}

Hence we may suppose $p+q\geq 1$. In this case, each $L_i$ is $2$-lef.

We prove the result by induction on $n$. Denote the first statement by $\HL_n$ and the second statement by $\HR_n$.

We first show that $\HR_{n-1} \Rightarrow \HL_n$.

Assume that $\phi \in H^{p, q} (X)$ satisfies
\begin{equation}\label{hl eq}
L_1 \cdot...\cdot L_{n-p-q} \cdot \phi =0.
\end{equation}
Under the assumption $\HR_{n-1}$, we need to verify that $\phi=0$.
To this end, we apply Proposition \ref{restriction prop}.
By Proposition \ref{restriction prop}, we can take a smooth hypersurface $V\in |mL_{n-p-q}|$ such that the restrictions ${L_1}_{|V},...,{L_{n-1-p-q}}_{|V}, {L_{n-p-q}}_{|V}$ satisfy the positivity condition for $\HR_{n-1}$ on $V$, that is,
\begin{itemize}
  \item for any $I \subset [n-1-p-q]$, ${L_I}_{|V}$ is $|I|+p+q$ lef on $V$;
  \item ${L_{n-p-q}}_{|V}$ is $p+q$ lef on $V$.
\end{itemize}

By restricting (\ref{hl eq}) to $V$, we get
\begin{equation}\label{hl restr}
{L_1}_{|V}\cdot...\cdot{L_{n-1-p-q}}_{|V}\cdot {L_{n-p-q}}_{|V}\cdot \phi_{|V} =0.
\end{equation}
Using the positivity condition on $V$ and $\HR_{n-1}$, we get that
\begin{equation}\label{hr n-1}
  c_{p,q} {L_1}_{|V}\cdot...\cdot{L_{n-1-p-q}}_{|V}\cdot \phi_{|V} \cdot \overline{\phi}_{|V} \geq 0
\end{equation}
with equality holds if and only if $\phi_{|V} =0$.

By the definition of $V$ and (\ref{hl eq}), it is clear that (\ref{hr n-1}) is an equality. Therefore, $\phi_{|V} =0$.
This implies $\phi =0$ since the restriction map
\begin{equation*}
  H^{p+q}(X, \mathbb{C}) \rightarrow H^{p+q} (V, \mathbb{C})
\end{equation*}
is injective by Lemma \ref{weakLef}.

Next, we show that $\HL_{n} \Rightarrow \HR_n$.

Recall that the primitive space is defined by
\begin{equation*}
  \prim^{p,q} (X) =\ker \{\Omega\cdot M : H^{p,q} (X, \mathbb{C})\rightarrow H^{n-q+1,n-p+1} (X, \mathbb{C})\}.
\end{equation*}
We claim that $H^{p,q}(X, \mathbb{C})$ admits a $Q$-orthogonal decomposition as follows:
\begin{equation}\label{or decom}
  H^{p,q}(X, \mathbb{C})= \prim^{p,q} (X) \oplus M\cdot H^{p-1,q-1} (X, \mathbb{C})
\end{equation}
with the convention that $H^{p-1,q-1} (X, \mathbb{C})=\{0\}$ when $p=0$ or $q=0$. 
To this end, we note that by the positivity assumption and Lemma \ref{sum lef max},
\begin{itemize}
  \item for any $I\subset [n-p-q]$, $L_I$ is $|I|+p+q$ lef, which is automatically $|I|+p+q-2$ lef.
  \item for any $I_1 \subset [n-p-q]$, $L_I := L_{I_1}+M$ is $\max\{|I_1| +p+q, p+q\}$ lef, which is automatically $|I|+p+q-2$ lef.
  \item for any $I_2 \subset [n-p-q]$, $L_I := L_{I_2}+M+M$ is $\max\{|I_2| +p+q, p+q\}$ lef, which is also $|I|+p+q-2$ lef.
\end{itemize}
Therefore, by $\HL_n$, the class $L_1\cdot...\cdot L_{n-p-q}\cdot M^2$ has Hard Lefschetz property, i.e.,
\begin{equation*}
  L_1\cdot...\cdot L_{n-p-q}\cdot M^2: H^{p-1,q-1} (X, \mathbb{C}) \rightarrow H^{n-q+1,n-p+1} (X, \mathbb{C})
\end{equation*}
is an isomorphism. In particular,
\begin{equation*}
  L_1\cdot...\cdot L_{n-p-q}\cdot M: M\cdot H^{p-1,q-1} (X, \mathbb{C}) \rightarrow H^{n-q+1,n-p+1} (X, \mathbb{C})
\end{equation*}
is an isomorphism. Combining the definition of $\prim^{p,q}$, we complete the proof of the claimed $Q$- decomposition for $  H^{p,q}(X, \mathbb{C})$.

As a consequence of the decomposition, we obtain that
\begin{equation}\label{dim prim}
  \dim \prim^{p,q}(X)=h^{p,q} -h^{p-1,q-1}.
\end{equation}

For $t\geq0$ and a fixed ample line bundle $A$, consider $L_1 +tA,...,L_{n-p-q} +tA, M+tA$ and the corresponding $\prim_t^{p,q}, Q_t$. Then by the mixed Hodge-Riemann bilinear relation for K\"ahler classes \cite{DN06, cattanimixedHRR}, we have
\begin{itemize}
  \item $\dim \prim_t ^{p,q}(X)=h^{p,q} -h^{p-1,q-1}=\dim \prim^{p,q}(X)$ for any $t>0$;
  \item $Q_t$ is positive definite on $\prim_t^{p,q}$ for any $t>0$.
\end{itemize}
By $\HL_n$,  $Q_t$ is non-degenerate on $H^{p,q}(X)$ for any $t\geq 0$. Therefore, $Q=\lim_{t\rightarrow 0} Q_t$ is positive definite on $\prim^{p,q}(X)$.

This finishes proof of the theorem.

\end{proof}

\bibliography{reference}

\providecommand{\bysame}{\leavevmode\hbox to3em{\hrulefill}\thinspace}
\providecommand{\MR}{\relax\ifhmode\unskip\space\fi MR }
\providecommand{\MRhref}[2]{%
  \href{http://www.ams.org/mathscinet-getitem?mr=#1}{#2}
}
\providecommand{\href}[2]{#2}
\begin{thebibliography}{{Sta}23}

\bibitem[ACG11]{GeometryOfAlgebraicCurvesVol2}
Enrico Arbarello, Maurizio Cornalba, and Phillip~A. Griffiths, \emph{Geometry
  of algebraic curves. {V}olume {II} with a contribution by {J}oseph {D}aniel
  {H}arris}, Grundlehren der mathematischen Wissenschaften [Fundamental
  Principles of Mathematical Sciences], vol. 268, Springer-Verlag, Heidelberg,
  2011. \MR{2807457}

\bibitem[Adi18]{adiprasito2019combinatorial}
Karim Adiprasito, \emph{Combinatorial {L}efschetz theorems beyond positivity},
  arXiv:1812.10454 (2018).

\bibitem[Cat08]{cattanimixedHRR}
Eduardo Cattani, \emph{Mixed {L}efschetz theorems and {H}odge-{R}iemann
  bilinear relations}, Int. Math. Res. Not. IMRN (2008), no.~10, Art. ID
  rnn025, 20. \MR{2429243}

\bibitem[dCM02]{decataldoLefsemismall}
Mark Andrea~A. de~Cataldo and Luca Migliorini, \emph{The hard {L}efschetz
  theorem and the topology of semismall maps}, Ann. Sci. \'{E}cole Norm. Sup.
  (4) \textbf{35} (2002), no.~5, 759--772. \MR{1951443}

\bibitem[dCM05]{catalMigHodgeMaps}
\bysame, \emph{The {H}odge theory of algebraic maps}, Ann. Sci. \'{E}cole Norm.
  Sup. (4) \textbf{38} (2005), no.~5, 693--750. \MR{2195257}

\bibitem[DN06]{DN06}
Tien-Cuong Dinh and Vi{\^e}t-Anh Nguy{\^e}n, \emph{The mixed {H}odge-{R}iemann
  bilinear relations for compact {K}\"ahler manifolds}, Geom. Funct. Anal.
  \textbf{16} (2006), no.~4, 838--849.

\bibitem[EV89]{viehwegvanishingnonvan}
H\'el\`ene Esnault and Eckart Viehweg, \emph{Vanishing and non-vanishing
  theorems}, Th\'eorie de Hodge - Luminy, Juin 1987 (Barlet D., Esnault H.,
  Elzein F., Verdier Jean-Louis, and Viehweg E., eds.), Ast\'erisque, no.
  179-180, Soci\'et\'e math\'ematique de France, 1989 (en). \MR{1042803}

\bibitem[GM88]{goreskymacphersonMorse}
Mark Goresky and Robert MacPherson, \emph{Stratified {M}orse theory},
  Ergebnisse der Mathematik und ihrer Grenzgebiete (3) [Results in Mathematics
  and Related Areas (3)], vol.~14, Springer-Verlag, Berlin, 1988. \MR{932724}

\bibitem[Gro65]{EGAiv}
Alexander Grothendieck, \emph{{\'{E}}l\'{e}ments de g\'{e}om\'{e}trie
  alg\'{e}brique. {IV}. {\'{e}}tude locale des sch\'{e}mas et des morphismes de
  sch\'{e}mas. {II}}, Inst. Hautes \'{E}tudes Sci. Publ. Math. \textbf{24}
  (1965), 231 pp. \MR{0199181}

\bibitem[Har77]{HartsAG}
Robin Hartshorne, \emph{Algebraic geometry}, Graduate Texts in Mathematics,
  vol.~52, Springer-Verlag, New York-Heidelberg, 1977. \MR{0463157}

\bibitem[Hir68]{hironakaprinciple}
Heisuke Hironaka, \emph{Smoothing of algebraic cycles of small dimensions},
  Amer. J. Math. \textbf{90} (1968), 1--54. \MR{224611}

\bibitem[Huh18]{huhICM2018}
June Huh, \emph{Combinatorial applications of the {H}odge-{R}iemann relations},
  Proceedings of the {I}nternational {C}ongress of {M}athematicians---{R}io de
  {J}aneiro 2018. {V}ol. {IV}. {I}nvited lectures, World Sci. Publ.,
  Hackensack, NJ, 2018, pp.~3093--3111. \MR{3966524}

\bibitem[Huh22]{huh2022combinatorics}
\bysame, \emph{Combinatorics and {H}odge theory}, Proceedings of the
  {I}nternational {C}ongress of {M}athematicians, 2022.

\bibitem[HX22]{huxiaohardlef2022Arxiv}
Jiajun Hu and Jian Xiao, \emph{Hard {L}efschetz properties, complete
  intersections and numerical dimensions}, arXiv:2212.13548 (2022).

\bibitem[Laz04]{lazarsfeldPosI}
Robert Lazarsfeld, \emph{Positivity in algebraic geometry. {I}}, Ergebnisse der
  Mathematik und ihrer Grenzgebiete. 3. Folge. A Series of Modern Surveys in
  Mathematics [Results in Mathematics and Related Areas. 3rd Series. A Series
  of Modern Surveys in Mathematics], vol.~48, Springer-Verlag, Berlin, 2004,
  Classical setting: line bundles and linear series. \MR{2095471}

\bibitem[Mum66]{Mumford}
David Mumford, \emph{Lectures on curves on an algebraic surface. {W}ith a
  section by {G}. {M}. {B}ergman}, Annals of Mathematics Studies, vol.~59,
  Princeton, N.J., 1966. \MR{0209285}

\bibitem[PT11]{takayamaIntermKodai}
Gianluca Pacienza and Shigeharu Takayama, \emph{On volumes along subvarieties
  of line bundles with nonnegative {K}odaira-{I}itaka dimension}, Michigan
  Math. J. \textbf{60} (2011), no.~1, 35--49. \MR{2785862}

\bibitem[RSW23]{ross2023duallylorent}
Julius Ross, Hendrik S\"{u}ss, and Thomas Wannerer, \emph{Dually {L}orentzian
  {P}olynomials}, arXiv:2304.08399 (2023).

\bibitem[Som78]{sommesMample}
Andrew~John Sommese, \emph{Submanifolds of {A}belian varieties}, Math. Ann.
  \textbf{233} (1978), no.~3, 229--256. \MR{466647}

\bibitem[SS85]{sommeseShiffmVanishingBook}
Bernard Shiffman and Andrew~John Sommese, \emph{Vanishing theorems on complex
  manifolds}, Progress in Mathematics, vol.~56, Birkh\"{a}user Boston, Inc.,
  Boston, MA, 1985. \MR{782484}

\bibitem[{Sta}23]{Stacks}
The {Stacks project authors}, \emph{The {S}tacks {P}roject}, 2023,
  \url{http://stacks.math.columbia.edu.}

\bibitem[Tim98]{timorinMixedHRR}
V.~A. Timorin, \emph{Mixed {H}odge-{R}iemann bilinear relations in a linear
  context}, Funktsional. Anal. i Prilozhen. \textbf{32} (1998), no.~4, 63--68,
  96. \MR{1678857}

\bibitem[Wil18]{williamECM}
Geordie Williamson, \emph{The {H}odge theory of the {H}ecke category}, European
  {C}ongress of {M}athematics, Eur. Math. Soc., Z\"{u}rich, 2018, pp.~663--683.
  \MR{3890447}

\bibitem[Xia21]{xiaomixedHRR}
Jian Xiao, \emph{Mixed {H}odge-{R}iemann bilinear relations and
  {$m$}-positivity}, Sci. China Math. \textbf{64} (2021), no.~7, 1703--1714.
  \MR{4280377}

\end{thebibliography}
\bibliographystyle{amsalpha}

\bigskip

\bigskip

\noindent
\textsc{Tsinghua University, Beijing 100084, China}\\
\noindent
\verb"Email: hujj22@mails.tsinghua.edu.cn"\\
\noindent
\verb"Email: sjshang@mail.tsinghua.edu.cn"\\
\noindent
\verb"Email: jianxiao@tsinghua.edu.cn"

\end{document}